\documentclass[10pt]{amsart}

\usepackage{amsmath, amsfonts, amsthm, amssymb, graphicx, fullpage, 
enumerate}

\newtheorem{theorem}{Theorem}      
\newtheorem{lemma}{Lemma}

\newtheorem{proposition}{Proposition}    
\newtheorem*{main theorem}{Main Theorem}   
\newtheorem*{thmA}{Theorem A}  
\newtheorem*{thmB}{Theorem B} 
\newtheorem*{thmC}{Theorem C}  
\theoremstyle{remark}
 
\theoremstyle{definition}
\newtheorem{definition}{Definition}

\def\N{\mathbb{N}}     
\def\Q{\mathbb{Q}}     
  
\def\Z{\mathbb{Z}}    
\def\C{\mathbb{C}} 
\def\P{\mathcal{P}}

\def\bar#1{\overline{#1}} 

\begin{document}
\title{Improved Bounds on S\'ark\"ozy's Theorem\\ for Quadratic Polynomials}
\author{Mariah Hamel\quad\quad\quad Neil Lyall \quad\quad\quad Alex Rice }

\address{D\'epartement de Math\'ematiques et de Statistique, Universit\'e de Montr\'eal, CP 6128, Centre-ville, Montr\'eal, QC H3C 3J7, Canada}
\email{mhamel@dms.umontreal.ca}
\address{Department of Mathematics, The University of Georgia, Athens, GA 30602, USA}
\email{lyall@math.uga.edu}
\address{Department of Mathematics, The University of Georgia, Athens, GA 30602, USA}
\email{arice@math.uga.edu} 


\begin{abstract} We extend the best known bound on the largest subset of $\{1,2,\dots,N\}$ with no square differences to the largest possible class of quadratic polynomials.  
\end{abstract}
\maketitle   
\setlength{\parskip}{5pt} 
\section{Introduction} 
\subsection{Background and previous results}
Lov\'asz conjectured that any set of natural numbers of positive upper density\footnote{ \ A set $A\subseteq \N$ is said to have positive upper density if
$\limsup_{N \to \infty} \frac{|A\cap[1,N]|}{N}>0,$ where $[1,N]$ denotes $\{1,2,\dots,N\}$.} necessarily contains two distinct elements which differ by a perfect square. This conjecture was proven in the late 1970s independently by S\'ark\"ozy and Furstenberg.
Furstenberg \cite{Furst} used ergodic theory and obtained a purely qualitative result, proving the conjecture exactly as stated above. S\'ark\"ozy, however, obtained a stronger, quantitative result by employing a Fourier analytic density increment strategy which utilized the Hardy-Littlewood circle method and was inspired by Roth's proof of the analogous conjecture for three-term arithmetic progressions \cite{Roth}.

\begin{thmA}[S\'ark\"ozy, \cite{sarkozy}]  If $A\subseteq [1,N]$ and $n^2 \notin A-A$ for all $n\in \N$, then\begin{equation} \label{sarkb} \frac{|A|}{N} \ll \Big(\frac{(\log \log N)^2}{\log N}\Big)^{1/3}. \end{equation}
\end{thmA}  
 
In this and the following theorems, $A-A$ denotes the difference set $\{a-a':a,a' \in A\}$, the symbol $\ll$ denotes ``less than a constant times'', and we implicitly assume that $N$ is large enough to make the right hand side of the inequalities defined and positive. The best known quantitative bound for the density of a subset $A\subseteq[1,N]$ with no square differences was obtained by Pintz, Steiger, and Szemer\'edi \cite{PSS}.

\begin{thmB}[Pintz, Steiger, and Szemer\'edi, \cite{PSS}]
If $A\subseteq [1,N]$ and $n^2 \notin A-A$ for all $n\in \N$, then\begin{equation} \label{pssb} \frac{|A|}{N} \ll (\log N)^{-c\log \log \log \log N}, \end{equation}
with $c=1/12$.
\end{thmB}  

A natural generalization of S\'ark\"ozy's theorem is the replacement of the squares with the image of a more general integer polynomial. In \cite{BPPS}, for example, Balog, Pelik\'an, Pintz, and Szemer\'edi modified the argument used in \cite{PSS} to obtain the same bounds with squares replaced by perfect $k^{\text{th}}$-powers for an arbitrary fixed $k\in \N$. In fact, they improved the constant $c$ in the exponent from $1/12$ to $1/4$. 
 
However, it is not the case that an analogous result can be obtained for an arbitrary polynomial, even in a qualitative sense. Given a polynomial $f\in \Z[x]$, it is clearly necessary that $f$ has a root modulo $q$ for every $q \in \N$, as otherwise there would be a set of the form $q\N$ of positive density with no differences in the image of $f$. It follows from a theorem of Kamae and Mend\`es France \cite{KMF} that this condition is also sufficient, and in this case we say that $f$ is an \textit{intersective polynomial}. Examples of intersective polynomials include any polynomial with an integer root and any polynomial with two rational roots with coprime denominators. However, there exist intersective polynomials with no rational roots, for example $(x^3-19)(x^2+x+1)$.

The first broad quantitative generalization of Theorem A beyond monomials was obtained  by Slijep\v{c}evi\'c \cite{Slip}, who showed triple logarithmic decay in the case of polynomials with an integer root. Lyall and Magyar \cite{LM}  obtained a stronger, single logarithmic bound in the integer root case as a corollary of a higher dimensional result (see also \cite{LM2}).  

The best bounds for an arbitrary intersective polynomial are due to Lucier, who successfully adapted the density increment procedure by utilizing $p$-adic roots and allowing the polynomial to change at each step of the iteration.

\begin{thmC}[Lucier, \cite{Lucier}] Suppose $f\in \Z[x]$ is an intersective polynomial of degree $k$. If $A \subseteq [1,N]$ and $f(n) \notin A-A$ for all $n \in \N$ with $f(n) \neq 0$, then \begin{equation*} \frac{|A|}{N} \ll \Big(\frac{(\log\log N)^{\mu}}{\log N}\Big)^{1/(k-1)}, \quad  \mu=\begin{cases}3 &\text{if }k=2 \\ 2 &\text{if }k>2 \end{cases}, \end{equation*} where the implied constant depends only on $f$.
\end{thmC}

\subsection{Main result of this paper}
In this paper, we combine Lucier's modified density increment strategy with the methods of \cite{PSS} and \cite{BPPS} in the special case of $k=2$. We also improve the constant in the exponent from $1/4$ to $1/\log 3$, the natural limit of the method as remarked in \cite{BPPS}, obtaining the following result. 
 
\begin{theorem} \label{main} Suppose $f \in \Z[x]$ is an intersective quadratic polynomial. If $A \subseteq [1,N]$ and $f(n) \notin A-A$ for all $n \in \N$ with $f(n) \neq 0$, then \begin{equation}\label{conc} \frac{|A|}{N} \ll (\log N)^{-\rho\log\log\log\log N} \end{equation}for any $\rho<1/\log 3$, where the implied constant depends only on $f$ and $\rho$.
\end{theorem}  

It is a pleasing consequence of Theorem \ref{main} and the previous results of \cite{PSS} and \cite{BPPS} that the primes, and even rather sparse subsets thereof, contain the desired arithmetic structure for any monomial or intersective quadratic based on density considerations alone. While the $1/\log N$ density barrier has not been broken for an arbitrary intersective polynomial, recent work of L\^{e} \cite{Le} uses Lucier's ideas together with Green's transference principle to show that for any intersective polynomial $f$, a subset of the primes of positive relative upper density is guaranteed to contain two distinct elements whose difference lies in the image of $f$.

\subsection{Remark on intersective quadratic polynomials}
It is worth pointing out that while the intersective condition can be somewhat mysterious and difficult to check for a general polynomial, this is not the case when restricted to degree $2$. 

\begin{proposition} \label{intquad} A quadratic polynomial $f \in \Z[x]$ is intersective if and only if $f$ has rational roots with coprime denominators. In other words, \begin{equation*} f(x)=a(\alpha x + \beta)(\gamma x + \lambda), \text{ } a,\alpha,\beta,\gamma,\lambda \in \Z, \text{ }(\alpha,\beta)=(\gamma,\lambda)=(\alpha,\gamma)=1. \end{equation*} 
\end{proposition} 

\noindent While it follows from Theorem 1 of \cite{BB} that an intersective polynomial with no rational roots must have degree at least five, Proposition \ref{intquad} can be directly shown more elementarily by noting that a polynomial is intersective if and only if it has a root in the $p$-adic integers for every prime $p$, then applying the quadratic formula over an appropriate field of $p$-adic numbers, and we include a short proof in Appendix \ref{B1}. While this characterization is not essential to the argument, it will allow us to greatly simplify some of Lucier's work. For example, we will avoid further discussion of $p$-adic numbers altogether.

\begin{center} \textbf{Acknowledgement} \end{center}

The authors would like to acknowledge Julia Wolf, whose exposition in \cite{julia} of the sensitive and initially intimidating argument in \cite{PSS} we found most helpful.

\section{Preliminary Notation: The Fourier Transform and the Circle Method}
We identify  subsets of the interval $[1,N]$ with subsets of the finite group $\Z_N=\Z/N\Z$, on which we utilize the normalized discrete Fourier transform. Specifically, for a function $F: \Z_N \to \C$, we define $\widehat{F}: \Z_N \to \C$ by \begin{equation*} \widehat{F}(t) = \frac{1}{N} \sum_{x \in \Z_N} F(x)e^{-2 \pi ixt/N}. \end{equation*} 

We analyze the behavior of the Fourier transform using the Hardy-Littlewood circle method, decomposing the nonzero frequencies  into two pieces: the roots of unity which are close to rationals with small denominator, and those which are not.

\begin{definition}\label{arcs} Given $N\in \N$ and a parameter $K>0$, we define, for each $q\in \N$ and $a \in [1,q]$, 
\begin{equation*} \mathbf{M}_{a,q}(K)=  \left\{t \in \Z_N: \left|\frac{t}{N}-\frac{a}{q}\right| < \frac{K}{N} \right\} \text{ \ and \ } \mathbf{M}_q(K) = \bigcup_{(a,q)=1}\mathbf{M}_{a,q}(K)\setminus \{0\}, \end{equation*} 
where the absolute value is on the circle, i.e. $0$ and $1$ are identified. We then define $\mathfrak{M}$, the \emph{major arcs}, by \begin{equation*} \mathfrak{M}(K) = \bigcup_{q=1}^{K^2} \mathbf{M}_q(K), \end{equation*} and $\mathfrak{m}(K)$, the \emph{minor arcs}, by $\mathfrak{m}(K) = \Z_N \setminus (\mathfrak{M}(K) \cup \{0\})$. It is important to note that as long as $2K^5<N$, we have that $a/q \neq b/r$ implies $\mathbf{M}_{a,q} \cap \mathbf{M}_{b,r} = \emptyset$ whenever $q,r \leq K^2$.
\end{definition} 

\noindent \textbf{Remark on notation.} We note that  the objects defined above certainly depend on $N$, despite its absence from the notation. In practice, $N$ should always be replaced with the size of the appropriate ambient group. When considering a set $A$, we let $A(x)$ denote the characteristic function of $A$, and we use the letters $C$ and $c$ to denote appropriately large or small constants which can change from line to line. 

\section{Overview of the argument}  
The underlying philosophy of this and many related results is that certain types of non-random phenomena in a set of integers should be detectable in the Fourier transform of the characteristic function of the set. That information about the transform can then be used to obtain some more explicit structural information about the set, such as increased density on a long arithmetic progression, and eventually provide an upper bound on its size. 

More specifically, we define $$I(f)=\{f(n)>0: n\in \N\}$$ for a polynomial $f\in\Z[x]$ with positive leading term. If $(A-A)\cap I(f)=\emptyset$ for a set $A\subseteq [1,N]$, one can apply the circle method and Weyl sum estimates to show that this unexpected behavior implies substantial $L^2$  mass of $\widehat{A}$ over nonzero frequencies near rationals with small denominator. At this point, there are multiple paths  to take in order to obtain the desired structural information. 

The original method of S\'ark\"ozy \cite{sarkozy}, as well as that of Lucier \cite{Lucier} and Lyall and Magyar \cite{LM}, is to use the pigeonhole principle to conclude that there is one single denominator $q$ such that $\widehat{A}$ has $L^2$ concentration around  rationals with denominator $q$. From this information, one can conclude that $A$ has increased density on a long arithmetic progression with step size an appropriate multiple of $q$, for example $q^2$ in the classical case. 
By translating and scaling the intersection of $A$ with this progression, one obtains a new subset $A'$ of a slightly smaller interval with significantly greater density. In addition, if $f$ is an intersective polynomial, $A'$ inherits non-random behavior from the fact that  $(A-A)\cap I(f)=\emptyset$. In the case that $f$ is a monomial, $A'$ actually inherits the identical property, but more generally one sees that $(A'-A')\cap I(h) = \emptyset$ for a potentially different intersective polynomial $h$ obtained from $f$. One then shows that if the density of the original set $A$ was too large, then this process could be iterated enough times for the density to surpass $1$, hence obtaining a contradiction. 

Pintz, Steiger, and Szemer\'edi \cite{PSS}  observed that pigeonholing to obtain a single denominator $q$ is a potentially wasteful step. We follow their approach, observing the following dichotomy: 
\begin{quotation}
\noindent \textbf{Case 1.} There is a single denominator $q$ such that $\widehat{A}$ has extremely high $L^2$ concentration, greater than yielded by the pigeonhole principle, around rationals with denominator $q$. This leads to a large density increment on a long arithmetic progression. 
\medskip

\noindent \textbf{Case 2.} The $L^2$ mass of $\widehat{A}$ on the major arcs is spread over many denominators. In this case, an iteration procedure using the ``combinatorics of rational numbers" can be employed to build a large collection of frequencies at which $\widehat{A}$ is large, then Plancherel's identity is applied to bound the density of $A$.
\end{quotation}

Philosophically, Case 1 provides more structural information about the original set $A$ than Case 2 does. The downside is that the density increment procedure yields a new set and potentially  a new polynomial, while the iteration in Case 2 leaves these objects fixed. 
With these cases in mind, we can now outline the argument, separated into two distinct phases. 
\begin{quotation}
\noindent \textbf{Phase 1} (The Outer Iteration): Given a set $A$ and an intersective quadratic polynomial $f$ with $(A-A)\cap I(f)=\emptyset,$ we ask if the set falls into Case 1 or Case 2 described above. If it falls into Case 2, then we proceed to Phase 2. 
\medskip

\noindent
If it falls into Case 1, then the density increment procedure yields a new subset $A_1$ of a slightly smaller interval with significantly greater density, and an intersective quadratic $f_1$ with slightly larger coefficients and $(A_1-A_1) \cap I(f_1) = \emptyset$. 
We can then iterate this process as long as the resulting interval is not too small, and the dichotomy holds as long as the coefficients of the corresponding polynomial are not too large.
\medskip

\noindent
One can show that if the resulting sets remain in Case 1, and the process iterates until the interval shrinks down or the coefficients grow to the limit, then the density of the original set $A$ must have satisfied a bound stronger than the one purported in Theorem \ref{main}.
\medskip

\noindent
Contrapositively, we assume that the original density does not satisfy this stricter bound, and we conclude that one of the sets yielded by the density increment procedure must lie in a large interval, have no differences in the image of a polynomial with small coefficients, and fall into Case 2. We call that set $B\subseteq[1,L]$ and the corresponding polynomial $h$.
\medskip
\end{quotation}

We now have a set $B \subseteq [1,L]$ and a quadratic polynomial $h$ with $(B-B)\cap I(h)= \emptyset$ which fall into Case 2, so we can adapt the strategy of \cite{PSS} and \cite{BPPS}. 

\begin{quotation} \textbf{Phase 2} (The Inner Iteration): We prove that given a frequency $s\in \Z_L$ with $s/L$ close to a rational $a/q$ such that $\widehat{B}(s)$ is large, there are lots of nonzero frequencies $t\in \Z_L$ with $t/L$ close to rationals $b/r$ such that $\widehat{B}(s+t)$ is almost as large. This intuitively indicates that a set $P$ of frequencies associated with large Fourier coefficients can be blown up to a much larger set $P'$ of frequencies associated with nearly as large Fourier coefficients. 
\medskip

\noindent
The only obstruction to this intuition is the possibility that there are many pairs $(a/q,b/r)$ and $(a'/q',b'/r')$ with $a/q+b/r = a'/q'+b'/r'$. Observations made in \cite{PSS} and \cite{BPPS} on the combinatorics of rational numbers demonstrate that this potentially harmful phenomenon can not occur terribly often. 
\medskip

\noindent
Starting with the trivially large Fourier coefficient at $0$, this process is applied as long as certain parameters are not too large, and the number of iterations is ultimately limited by the growth of the divisor function. Once the iteration is exhausted, we use the resulting set of large Fourier coefficients and Plancherel's Identity to get the upper bound on the density of $B$, which is by construction larger than the density of the original set $A$, claimed in Theorem \ref{main}.
\end{quotation}

\section{Reduction of Theorem \ref{main} to Two Lemmas}
To make the strategy outlined in Section 3 precise, we fix an intersective quadratic $f$ and a set $A \subseteq [1,N]$ with $|A|=\delta N$ and $f(n) \notin A-A$ for all $n\in \N$ with $f(n)\neq 0$. By the symmetry of difference sets, we can assume without loss of generality that $f$ has positive leading term, and we see in particular that $(A-A)\cap I(f)=\emptyset$. 

We also fix an arbitrary $\epsilon > 0$ and set $Q=(\log N)^{\epsilon \log\log\log N}$, and we will prove Theorem \ref{main} with $\rho=(1-11\epsilon)/\log 3$. At any point we are free to insist that  $N$ is sufficiently large with respect to $f$ and $\epsilon$, as this will only affect the implied constant in (\ref{conc}),  so for convenience we will take these to be  perpetual implicit hypotheses and abstain from including them further.  From this point on, we will allow all of our constants to depend on $f$ and $\epsilon$. 

\subsection{Two Key Lemmas}
We first reduce Theorem \ref{main} to two key lemmas, the first of which corresponds to Phase 1 outlined above, and the second of which corresponds to Phase 2. 
 
\begin{lemma} \label{outer} If \begin{equation} \label{behrend} \delta \geq e^{-(\log N)^{\epsilon/8}},\end{equation} then there exists  $B \subseteq [1,L]$ satisfying $L \geq N^{.99}$, $|B|/L=\sigma \geq \delta$, and $(B-B) \cap I(h) = \emptyset$, where \begin{equation*} h(x)= \begin{cases}ax^2, &\text{if } f(x)=a(x-b)^2\\ f(r+dx)/d, &\text{else} \end{cases},  \end{equation*}$d \leq N^{.01}$, and $r\in (-d,0]$ is a root of $f$ modulo $d$. Further, $B$ satisfies $|B\cap[1,L/2]|\geq \sigma L/3$ and  \begin{equation} \label{maxmass} \max_{q\leq Q}\sum_{t\in \mathbf{M}_q(Q)}  |\widehat{B}(t)|^2 \leq \sigma^2(\log N)^{-1+\epsilon}. \end{equation} 
\end{lemma}  

We note that $e^{-(\log N)^{\epsilon/8}} \ll (\log N)^{-\log\log\log\log  N}$. In particular, if hypothesis (\ref{behrend}) is not satisfied, then Theorem \ref{main} is already more than true. The next lemma corresponds to the iteration scheme in which a set of large Fourier coefficients from distinct major arcs is blown up in such a way that the relative growth of the size of the set is much greater than the relative loss of pointwise mass.

\begin{lemma} \label{itn} Suppose $B\subseteq [1,L]$ and $h\in \Z[x]$ are as in the conclusion of Lemma \ref{outer}, let $B_1=B\cap[1,L/2]$, and suppose $\sigma \geq Q^{-1/6}$. Given $U,V,K \in \N$ with $\max \{U,V,K\}\leq Q^{1/6}$ and a set \begin{equation*} P \subseteq \left\{t \in \bigcup_{q=1}^{V} \mathbf{M}_q(K)\cup\{0\} : |\widehat{B_1}(t)|  \geq \frac{\sigma}{U} \right\} \end{equation*} satisfying \begin{equation} \label{dis}| P \cap \mathbf{M}_{a,q}(K)| \leq 1 \text{ \ whenever \ } q \leq V,\end{equation} there exist $U',V',K' \in \N$ with $\max\{U',V',K'\} \ll (\max \{U,V,K\})^3\sigma^{-5/2}$ and a set  \begin{equation} \label{P'1} P' \subseteq \left\{t \in \bigcup_{q=1}^{V'} \mathbf{M}_q(K')\cup\{0\} : |\widehat{B_1}(t)|  \geq \frac{\sigma}{U'} \right\}\end{equation} satisfying \begin{equation} \label{P'2} | P' \cap \mathbf{M}_{a,q}(K')| \leq 1 \text{ \ whenever \ } q\leq V'\end{equation} and \begin{equation} \label{P'3} \frac{|P'|}{(U')^2} \geq \frac{|P|}{U^2}(\log N)^{1-10\epsilon}.\end{equation}
\end{lemma}

\subsection{Proof that Lemmas \ref{outer} and \ref{itn} imply Theorem \ref{main}} In order to prove Theorem \ref{main}, we can assume that \begin{equation*}\delta \geq (\log N)^{-\log\log\log\log N}. \end{equation*} Therefore, Lemma \ref{outer} produces a set $B$ of density $\sigma \geq \delta$ with the stipulated properties, and we set $P_0=\{0\}$, $U_0=3$, and $V_0=K_0=1$. Lemma \ref{itn} then yields, for each $n$, a set $P_n$ with parameters $U_n,V_n,K_n$ such that \begin{equation*}\max\{U_n,V_n,K_n\} \leq (\log N)^{3^{n+1}\log\log\log\log N} \end{equation*} and \begin{equation*}\frac{1}{\sigma} \geq \frac{1}{\sigma^2}\sum_{t\in P_n} |\widehat{B_1}(t)|^2 \geq \frac{|P_n|}{U_n^2} \gg (\log N)^{n(1-10\epsilon)}, \end{equation*} where the left-hand inequality comes from Plancherel's Identity, as long as $\max\{U_n,V_n,K_n\} \leq Q^{1/6}$. This holds with $n=(1-\epsilon)(\log\log\log\log N) /\log 3$, as $3^{n+1} \leq (\log\log\log N)^{1-\epsilon/2}$, and Theorem \ref{main} follows. 
\qed

\section{The Outer Iteration: Proof of Lemma \ref{outer}}

Recall that we have fixed $A \subseteq [1,N]$, an intersective quadratic $f\in \Z[x]$, and $\epsilon>0$, and in this section we also assume the bound (\ref{behrend}) on the density $\delta$. As previously mentioned, we will apply the modified density increment strategy described in \cite{Lucier}, which allows for the polynomial to change at each stage of the iteration. The following definition completely describes all of the polynomials that we could potentially encounter.
 
\begin{definition}\label{aux} For each $d\in\N$, we fix an integer $r_d \in (-d,0]$ such that $f(r_d) \equiv 0$ mod $d$ and $r_d \equiv r_s$ mod $s$ whenever $s \mid d$, and we define the \emph{auxiliary polynomials} $f_d\in \Z[x]$ by \begin{equation*} f_d(x)= \begin{cases}ax^2, &\text{if } f(x)=a(x-b)^2\\ f(r_d+dx)/d, &\text{else} \end{cases}.\end{equation*}
\end{definition}  

\noindent One can find a collection of roots with the property stipulated in Definition \ref{aux} in the following way. If $f(x)=a(\alpha x + \beta)(\gamma x + \lambda)$ with $(\alpha,\gamma)=1$, partition the primes into $\P=\P_1 \cup \P_2$, with $p \nmid \alpha$ for all $p \in \P_1$ and $p \nmid \gamma$ for all $p \in \P_2$.\\\\ For each $d\in \N$, write $d=p_1^{a_1}\cdots p_{k}^{a_{k}}s_1^{j_1}\cdots s_{\ell}^{j_{\ell}}$ with $p_i \in \P_1$ and $s_i \in \P_2$. By the Chinese Remainder Theorem, there is a unique integer $r_d \in (-d,0]$ such that $r_d \equiv -\beta \alpha^{-1} \mod p_i^{a_i}$ and $r_d \equiv -\lambda \gamma^{-1} \mod s_n^{j_n}$ for all $1\leq i \leq k$ and $1\leq n \leq \ell$, and we see that this choice of $r_d$ meets the purported condition. One can easily show from the characterization in Proposition \ref{intquad} and the construction of the roots $r_d$ that each of these auxiliary polynomials are themselves intersective quadratics.

We now invoke the usual density increment lemma which states that $L^2$-concentration of the Fourier transform leads to increased density on a long progression, with the added observation that if the difference set misses the image of a polynomial, then the difference set of the resulting subset of a smaller interval misses the image of an appropriate auxiliary polynomial. The particular statement below follows from Lemma 20 of \cite{Lucier}, while the additional observation is made in Lemma 31 of the same paper.
\begin{lemma} \label{inc} Suppose $B\subseteq [1,L]$ with $|B|=\sigma L$ and $(B-B)\cap I(f_d)=\emptyset$. If $L\geq Q^4$ and 
\begin{equation}\label{Bmass} \sum_{t\in \mathbf{M}_q(Q)} |\widehat{B}(t)|^2 \geq \sigma^2(\log N)^{-1+\epsilon},\end{equation}
 for some $q \leq Q$, then there exists $B'\subseteq [1,L']$ satisfying $L' \gg \sigma L/Q^4,$ $(B'-B')\cap I(f_{qd})=\emptyset,$ and
\begin{equation*} 
|B'|/L'\geq  \sigma(1+(\log N)^{-1+\epsilon}/8).\end{equation*}

\end{lemma} 
 
\subsection{Proof of Lemma \ref{outer}} 



Setting $A_0$=$A$, $N_0=N$, $\delta_0=\delta$, and $d_0=1$, we iteratively apply Lemma \ref{inc}. This yields, for each $k$, a set $A_k \subseteq [1,N_k]$ with $|A_k|=\delta_kN_k$ and $(A_k-A_k)\cap I(f_{d_k})=\emptyset$ satisfying
\begin{equation} \label{Ndeld} N_k \geq (c\delta/Q^4)^kN, \quad \delta_k \geq \delta_{k-1}(1+(\log N)^{-1+\epsilon}/8), \quad \text{and} \quad d_k \leq Q^k
\end{equation}
as long as $N_k \geq Q^4$ and either
\begin{equation}\label{newmass}\max_{q\leq Q}\sum_{t\in \mathbf{M}_q(Q)} |\widehat{A_k}(t)|^2\geq \delta_k^2(\log N)^{-1+\epsilon}
\end{equation}
or $|A_k\cap[1,N_k/2]|<\delta_kN_k/3$, as the latter condition implies $A_k$ has density at least $3\delta_k/2$ on the interval $(N_k/2,N_k]$. We see that by (\ref{behrend}) and (\ref{Ndeld}), the density $\delta_k$ will exceed 1 after \begin{equation*}16\log(\delta^{-1})(\log N)^{1-\epsilon} \leq (\log N)^{1-\epsilon/2}\end{equation*} steps, hence $N_k <Q^4$ or  (\ref{newmass}) fails and $|A_k\cap[1,N_k/2]|\geq \delta_kN_k/3$ for some 
\begin{equation} \label{kbound} k\leq (\log N)^{1-\epsilon/2}.\end{equation}
However, we see that (\ref{behrend}), (\ref{Ndeld}), and (\ref{kbound}) imply 
 \begin{equation*}N_k \geq N(e^{-(\log N)^{\epsilon/4}})^{(\log N)^{1-\epsilon/2}} \geq Ne^{-(\log N)^{1-\epsilon/4}} \geq N^{.99}, \end{equation*}
so we set $B=A_k$, $L=N_k$, $\sigma = \delta_k$, and $h=f_{d_k}$, and we see further that
\begin{equation*} d_k \leq Q^{(\log N)^{1-\epsilon/2}}\leq e^{(\log N)^{1-\epsilon/4}} \leq N^{.01}, \end{equation*} 
as required.
\qed

\section{The Inner Iteration: Proof of Lemma \ref{itn}}
In this section, we let $B\subseteq [1,L]$ and $h \in \Z[x]$ be as in the conclusion of Lemma \ref{outer}, we let $B_1=B\cap[1,L/2]$, and we assume $\sigma \geq Q^{-1/6}$. We set $j=\max \{n \in \N: h(n)\leq 0 \}+1,$ taking $j=1$ if $h$ is strictly positive on $\N$, and $$M=\min\{ n\in \N: h(n) \geq L/3 \}-1.$$ One can see from the definition of $h$, the lower bound on $L$, and the upper bound on $d$ that $j$ is uniformly bounded in terms of the original polynomial $f$, while $M \gg N^{.49}$. 

\subsection{Proof of Lemma \ref{itn}}
\noindent Suppose we have a set $P$ with parameters $U,V,K$ as specified in the hypotheses of Lemma \ref{itn}, and fix an element $s \in P$. Since $(B-B)\cap I(h)=\emptyset$, we see that there are no solutions to \begin{equation*}a-b\equiv h(y) \mod{L}, \quad a\in B, \  b \in B_1, \  j\leq y \leq M. \end{equation*}  
Combined with the orthogonality of the characters, this implies
\begin{equation*} \sum_{t \in \Z_{L}} \widehat{B}(t)\bar{\widehat{B_1}(s+t)}S(t) =\frac{1}{LM^2}\sum_{\substack{x\in\Z_{L} \\ j\leq y \leq M}} yB(x+h(y))B_1(x)e^{2\pi i xs/L}= 0,\end{equation*} 
where $$S(t)=\frac{1}{M^2}\sum_{x=j}^M xe^{2\pi i h(x) t/L},$$which immediately yields 
\begin{equation} \label{mass} \sum_{t\in \Z_{L} \setminus \{0\}} |\widehat{B}(t)||\widehat{B_1}(s+t)||S(t)| \geq \Big|\sum_{t\in \Z_{L} \setminus \{0\}} \widehat{B}(t)\bar{\widehat{B_1}(s+t)}S(t)\Big| = \widehat{B}(0)\widehat{B_1}(s)S(0)\geq \sigma^2 /4U,\end{equation} 
since $|\widehat{B_1}(s)| \geq \sigma/U$ and $S(0)\geq 1/4.$ Letting $\eta=c_0\sigma/U$ for a constant $c_0>0$, it follows from traditional Weyl sum estimates and Lemmas 11 and 28 of \cite{Lucier} that 
\begin{equation}\label{Wmin} |S(t)| \leq C\eta  \leq \sigma/8U \quad \text{for all } t \in \mathfrak{m(\eta^{-1})}, \end{equation} 
provided $c_0 \leq 1/8C$, and \begin{equation} \label{Wmaj} |S(t)| \ll q^{-1/2}\min\{1, (L|t/L-a/q|)^{-1} \} \end{equation} 
if $t \in \mathbf{M}_{a/q}(\eta^{-1}),$ $(a,q)=1$, and $q\leq \eta^{-2}$. We discuss estimates (\ref{Wmin}) and (\ref{Wmaj}) in more detail in Appendix \ref{A}.

\noindent We have by (\ref{Wmin}), Cauchy-Schwarz, and Plancherel's Identity that 
\begin{equation*} \sum_{t \in \mathfrak{m}(\eta^{-1}) } |\widehat{B}(t)||\widehat{B_1}(t)||S(t)|  \leq \sigma^2/8U,  \end{equation*} 
which together with (\ref{mass}) yields 
\begin{equation} \label{Mmass} \sum_{t \in \mathfrak{M}(\eta^{-1})  } |\widehat{B}(t)||\widehat{B_1}(t)||S(t)|  \geq \sigma^2 /8U.  \end{equation} 
We now wish to assert that we can ignore those frequencies in the major arcs at which the transform of $B$ or $B_1$ is particularly small. In order to make this precise, we first need to invoke a weighted version of a well-known estimate on the higher moments of Weyl sums. Specifically, we have that  
\begin{equation} \label{moment} \sum_{t \in \Z_{L}} |S(t)|^6 \leq C, \end{equation} 
which can be seen by adapting the proof of Proposition 3.3 of \cite{LM2}, and we provide a proof in Appendix \ref{A}. Choosing a constant $0<c_1<(32C^{1/6})^{-3}$, where $C$ comes from (\ref{moment}), we define 
\begin{equation} \label{XY} X = \left\{ t\in \mathfrak{M}(\eta^{-1})  :\min\Big\{|\widehat{B}(t)|,|\widehat{B_1}(s+t)|\Big\} \leq c_1\sigma^{7/2}/U^3\right\} \text{ \ and \ } Y=\mathfrak{M}(\eta^{-1})\setminus X. \end{equation} 
Using H\"{o}lder's Inequality to exploit the sixth moment estimate on $S$, followed by Plancherel's Identity, we see that 
\begin{align*}\sum_{t \in X}|\widehat{B}(t)||\widehat{B_1}(s+t)||S(t)| &\leq \Big(\sum_{t \in X} |\widehat{B}(t)|^{6/5}|\widehat{B_1}(s+t)|^{6/5} \Big)^{5/6} \Big( \sum_{t \in \Z_{L}} |S(t)|^6 \Big)^{1/6}\\ &\leq \frac{c_1^{1/3}\sigma^{7/6}}{U}\Big(\sum_{t \in \Z_{L}} \max\left\{|\widehat{B}(t)|^2,|\widehat{B_1}(s+t)|^2\right\}\Big)^{5/6}\cdot C^{1/6} \\ & \leq \frac{\sigma^{7/6}}{32U}\Big(\sum_{t \in \Z_{L}} |\widehat{B}(t)|^2+|\widehat{B_1}(s+t)|^2\Big)^{5/6} \leq \sigma^2/16U,\end{align*}
and hence by (\ref{Mmass}) we have 
\begin{equation}\label{Ymass}  \sum_{t \in Y } |\widehat{B}(t)||\widehat{B_1}(s+t)||S(t)|  \geq \sigma^2 /16U. \end{equation} 
For $i,j,k \in \N$, we define 
\begin{equation*} \mathcal{R}_{i,j,k} = \left\{ a/q: (a,q)=1, \ 2^{i-1} \leq q \leq 2^i,  \text{ } \frac{\sigma}{2^j} \leq \max |\widehat{B}(t)| \leq \frac{\sigma}{2^{j-1}},\text{ } \frac{\sigma}{2^k} \leq \max |\widehat{B_1}(s+t)| \leq \frac{\sigma}{2^{k-1}} \right\}, \end{equation*}
where the maximums are taken over nonzero frequencies  $t\in \mathbf{M}_{a/q}(\eta^{-1})$. We see that we have 
\begin{equation}\label{rij}\sum_{a/q\in \mathcal{R}_{i,j,k}}\sum_{t\in \mathbf{M}_{a/q}(\eta^{-1})\setminus\{0\}} |\widehat{B}(t)||\widehat {B_1}(s+t)||S(t)| \ll |\mathcal{R}_{i,j,k}|\frac{\sigma^2}{2^{j}2^k}\max_{a/q \in \mathcal{R}_{i,j,k}}\sum_{t\in \mathbf{M}_{a/q}(\eta^{-1})} |S(t)|. \end{equation}
It follows from (\ref{Wmaj}) and the bound $U,\sigma^{-1}\leq Q^{1/6}$ that if $(a,q)=1$ and $q\leq \eta^{-2}$, then
\begin{equation*}\sum_{t\in \mathbf{M}_{a/q}(\eta^{-1})} |S(t)| \ll q^{-1/2}\log(Q),\end{equation*} 
hence by (\ref{rij}) we have 
\begin{equation} \label{rijs}\sum_{a/q\in \mathcal{R}_{i,j,k}}\sum_{t\in \mathbf{M}_{a/q}(\eta^{-1})\setminus\{0\}} |\widehat{B}(t)||\widehat{B_1}(s+t)||S(t)| \ll |\mathcal{R}_{i,j,k}|\frac{\sigma^2}{2^{j}2^k}2^{-i/2}\log(Q). \end{equation}  
By our definitions, the sets $\mathcal{R}_{i,j,k}$ exhaust $Y$ by taking $1\leq 2^i \leq \eta^{-2}$ and $1 \leq 2^j,2^k \leq U^3/c_1\sigma^{5/2}$, a total search space of size $\ll (\log Q)^3 $. Therefore, by (\ref{Ymass}) and (\ref{rijs}) there exist $i,j,k$ in the above range such that 
\begin{equation*}\frac{\sigma^2}{U(\log Q)^3} \ll|\mathcal{R}_{i,j,k}|\frac{\sigma^2}{2^j2^k}2^{-i/2}\log Q. \end{equation*}
In other words, we can set $V_s=2^{i}$, $W_s=2^j$, and $U_s=2^k$ and take an appropriate nonzero frequency from each of the pairwise disjoint major arcs specified by $\mathcal{R}_{i,j,k}$ to form a set 
\begin{equation*} P_s \subseteq \left\{t \in \bigcup_{q=V_s/2}^{V_s} \mathbf{M}_{q}(\eta^{-1}):  \text{ } |\widehat{B_1}(s+t)| \geq \frac{\sigma}{U_s} \right\} \end{equation*} 
which satisfies 
\begin{equation}\label{Ps} |P_s| \gg \frac{U_sW_sV_s^{1/2}}{U(\log Q)^4},  \quad |P_s\cap\mathbf{M}_{a,q}(\eta^{-1})|\leq 1 \text{ \ whenever \ } q \leq V_s,  \end{equation}  
and 
\begin{equation}\label{maxBW} \max_{ t \in \mathbf{M}_{a/q}(\eta^{-1})\setminus \{0\}} |\widehat{B}(t)| \geq \frac{\sigma}{W_s} \text{ whenever } q\leq V_s \text{ and }  \mathbf{M}_{a/q}(\eta^{-1})\cap P_s \neq \emptyset,
\end{equation} 
noting by disjointness that $a/q \in \mathcal{R}_{i,j,k}$  whenever $q\leq V_s$ and  $\mathbf{M}_{a/q}(\eta^{-1})\cap P_s \neq \emptyset $. 

\noindent We now observe that there is a subset $\tilde{P} \subseteq P$ with 
\begin{equation}\label{tilde}|\tilde{P}| \gg |P|/(\log Q)^3\end{equation}  
for which the triple $U_s,W_s,V_s$ is the same. We call  those common parameters $\tilde{U},\tilde{W}$ and $\tilde{V}$, respectively, and we can now foreshadow by asserting that the claimed parameters in the conclusion of Lemma \ref{itn} will be $U'=\tilde{U}$, $V'=\tilde{V}V$, and $K'=K+\eta^{-1}$, which do satisfy the purported bound. 

\noindent We let 
\begin{equation*} \mathcal{R} = \left\{\frac{a}{q}+\frac{b}{r}:  s\in \mathbf{M}_{a/q}(K) \text{ for some } s \in \tilde{P} \text{ and }  t \in \mathbf{M}_{b/r}(\eta^{-1}) \text{ for some } t \in P_s \right\}. \end{equation*}
By taking one frequency $s+t$ associated to each element in $\mathcal{R}$,  we form our set $P'$, which immediately satisfies conditions (\ref{P'1}) and (\ref{P'2}) from the conclusion of Lemma \ref{itn}. However, the crucial condition (\ref{P'3}) on $|P'|$, which by construction is equal to $|\mathcal{R}|$, remains to be shown. To this end, we invoke the work on the combinatorics of rational numbers found in \cite{PSS} and \cite{BPPS}.

\begin{lemma}[Lemma CR of \cite{BPPS}]\label{CR} \begin{equation*} |\mathcal{R}|\geq \frac{|\tilde{P}|(\min_{s\in \tilde{P}}|P_s|)^2}{\tilde{V}D\tau^8(1+\log V)},  \end{equation*} where \begin{equation*} D= \max_{r\leq \tilde{V}} \Bigl|\Bigl\{ b  : \ (b,r)=1, \  \mathbf{M}_{b/r}(\eta^{-1})\cap \bigcup_{s\in \tilde{P}} P_s\neq \emptyset \Bigr\}\Bigr|, \end{equation*} $\tau(q)$ is the divisor function and $\tau = \max_{q\leq V\tilde{V}} \tau(q)$.
\end{lemma}

\noindent It is a well-known fact of the divisor function that $\tau(n) \leq n^{1/\log\log n}$ for large $n$, and since $V\tilde{V} \leq Q$, we have that $\tau \leq (\log N)^{\epsilon}$. 
We also have from (\ref{maxmass}) that 
\begin{equation*} \sigma^2(\log N)^{-1+\epsilon} \geq \max_{r \leq Q} \sum_{t \in \mathbf{M}_r(Q)} |\widehat{B}(t)|^2 \geq \max_{r \leq \tilde{V}} \sum_{t \in \mathbf{M}_r(\eta^{-1})} |\widehat{B}(t)|^2 \geq \frac{\sigma^2}{\tilde{W}^2}D,\end{equation*}
where the last inequality follows from (\ref{maxBW}),  and hence 
\begin{equation}\label{D} D \leq \tilde{W}^2(\log N)^{-1 + \epsilon}.\end{equation}

 \noindent Combining the estimates on $\tau$ and $D$ with (\ref{Ps}), (\ref{tilde}), and Lemma \ref{CR}, we have\begin{equation*}|P'| \gg \frac{|P|}{(\log Q)^3} \frac{\tilde{U}^2\tilde{W}^2\tilde{V}}{U^2(\log Q)^8} \frac{(\log N)^{1-\epsilon}}{\tilde{V}\tilde{W}^2(\log N)^{8\epsilon}(\log Q)} \geq \tilde{U}^2\frac{|P|}{U^2}(\log N)^{1-10\epsilon}.\end{equation*} \\
Recalling that we set $U'=\tilde{U}$, the lemma follows.
\qed

\appendix
\section{Exponential Sum Estimates: Proof of (\ref{Wmin}), (\ref{Wmaj}), and (\ref{moment})}\label{A}
Throughout Appendix A we write $h(x)=\alpha x^2+ \beta x + \gamma$. To begin, we invoke some Weyl sum estimates.
 
\begin{lemma}\label{major} If $t\in \Z_{L}$ and $t/L=a/q+\lambda$ with $q\leq M^{0.1}$, $(a,q)=1$, and $|\lambda|<M^{-1.9}$, then
\begin{equation*}   
S(t) = \frac{1}{qM^2}G(a,q)\int_1^M x e^{2 \pi i h(x) \lambda} dx+ O(M^{-0.7}),
\end{equation*}
where 
\begin{equation*}
G(a,q)=\sum_{r=0}^{q-1}e^{2 \pi i h(r)a/q}.
\end{equation*}
\end{lemma}
\begin{lemma} \label{minor}
If $t\in \Z_{L}$, $|t/L-a/q|<1/q^2$, and $(a,q)=1$, then
\begin{equation*}
|S(t)| \ll \log M (\alpha/q+\alpha/M+q/M^2)^{1/2}
\end{equation*}
\end{lemma}

Lemma \ref{major} is a weighted version of the traditional major arc asymptotic for Weyl sums, and in particular follows from Lemma 11 of \cite{Lucier}. Lemma \ref{minor} follows from the standard Weyl Inequality for quadratic polynomials (see \cite{Mont} for example) and summation by parts. 

The presence of the Gauss sum $G(a,q)$ in the conclusion of Lemma \ref{major} indicates that our estimates could be irreparably damaged if the coefficients of $h$ share large factors. The following observation of Lucier ensures that this feared scenario does not occur.

\begin{proposition}\label{content} If $f(x)=(\alpha x + \beta)(\gamma x + \lambda)$ with $\alpha,\beta,\gamma,\lambda \in \Z$ and $f$ does not have a double root, then for any $d\in \N$, \begin{equation*} \textnormal{cont}(f_d) \leq |\alpha\lambda-\beta\gamma|, \end{equation*} where $$\textnormal{cont}(a_0+a_1x+a_2x^2)=\gcd(a_1,a_2).$$ 
\end{proposition} 

\noindent We note that in the case $f(x)=a(x-b)^2$ excluded by the hypotheses, we trivially have cont$(f_d)=a$ for all $d$. While Proposition \ref{content}, which is a special case of Lemma 28 of \cite{Lucier}, is pleasingly precise, we will only use that $\textnormal{cont}(h)$ is uniformly bounded in terms of the original polynomial $f$.  Again, the degree $2$ case is considerably simpler than the analogous result for a general intersective polynomial which may or may not have rational roots, and we include an elementary proof of Proposition \ref{content} in Appendix \ref{B2}. 

We will also need some additional facts about the polynomial $h$ which follow from its construction as an auxiliary polynomial of $f$ in Lemma \ref{outer}. Specifically, the bounds on $d$ and $L$ in Lemma \ref{outer} and the definition of $M$ tell us  that $\max \{\alpha, |\beta|, |\gamma| \} \ll d \leq N^{.01},$ hence 
\begin{equation} \label{Mlow} M \gg N^{.49}
\end{equation}
and \begin{equation} \label{cofmax} \max \{\alpha, |\beta|, |\gamma| \} < M^{.03}. \end{equation}
Also, $ \alpha \gg |\beta|+|\gamma|$, and therefore 
\begin{equation} \label{aM2} \alpha M^2 \geq L/4.
\end{equation}
Finally, by (\ref{behrend}) and (\ref{Mlow}), we have that 
\begin{equation}\label{etaUp} \eta >M^{-.01}, \end{equation} and we are ready to establish estimates (\ref{Wmin}), (\ref{Wmaj}), and (\ref{moment}).

\subsection{Proof of  (\ref{Wmaj})} By (\ref{etaUp}), we see that for $U,\sigma^{-1}\leq Q^{1/6}$ the hypotheses of Lemma \ref{major} are comfortably satisfied whenever $t\in \mathbf{M}_{a/q}(\eta^{-1})$, $(a,q)=1$, and $q\leq \eta^{-2}$.  

\noindent For the Gauss sum $G(a,q)$, we use the well known estimate 
\begin{equation} \label{gauss}
|G(a,q)| \ll  (\text{cont}(h)q)^{1/2},
\end{equation}
which, for example, is a special case of Lemma 6 in \cite{Lucier}. Combining (\ref{gauss}) and (\ref{aM2}) with Lemma \ref{major} and Proposition \ref{content}, it suffices to show 
\begin{equation} \label{int}
\Big| \int_1^M xe^{2 \pi i h(x) \lambda} dx \Big| =\Big| \int_1^M xe^{2 \pi i(\alpha x^2+\beta x) \lambda} dx \Big| \leq \min \left\{M^2, (\alpha |\lambda|)^{-1}\right\}.
\end{equation}
The equality and the first of the two implicit inequalities are trivial, and in particular hold with $\lambda=0$. For $\lambda \neq 0$, we first ignore the linear term in the polynomial by observing 
\begin{equation*} \Big| \int_1^M xe^{2 \pi i(\alpha x^2+\beta x) \lambda} - xe^{2\pi i \alpha x^2 \lambda} dx \Big| \leq M  \int_1^M |e^{2\pi i \beta x \lambda}-1| dx \leq 2\pi M^3|\beta||\lambda| \leq M^{-0.7}/|\lambda|,
\end{equation*} since $|\beta|< M^{.03}$ by (\ref{cofmax})  and $|\lambda|<M^{-1.9}$ by assumption. 

\noindent For the main term, we change variables $(y:=\alpha x^2)$ to see 
\begin{equation*} \int_1^M  xe^{2\pi i \alpha x^2 \lambda} dx = \frac{1}{2\alpha} \int_{\alpha}^{\alpha M^2} e^{2\pi i y\lambda} dy = \frac{1}{2\alpha} \Big( \frac{e^{2\pi i \alpha M^2\lambda}-e^{2\pi i \alpha\lambda}}{2\pi i \lambda} \Big),
\end{equation*} which in absolute value is clearly at most $(2\pi \alpha|\lambda|)^{-1}$, and the estimate follows.
\qed

\subsection{Proof of (\ref{Wmin})} Fixing $t\in \mathfrak{m}(\eta^{-1})$, we have by the pigeonhole principle that there exist $1\leq q \leq M^{1.9}$ and $(a,q)=1$ with $|t/L-a/q| < 1/(qM^{1.9})$. 
If $\eta^{-2} \leq q \leq M^{0.1}$, then Lemma \ref{major} with the trivial bound on the integral, (\ref{gauss}), and Proposition \ref{content} immediately yield the desired estimate.
If $M^{0.1}\leq q \leq M^{1.9}$, then (\ref{cofmax}) and Lemma \ref{minor} imply 
\begin{equation*} |S(t)| \leq M^{-.03},
\end{equation*} which by (\ref{etaUp}) is stronger than the required estimate.
If $1\leq q \leq \eta^{-2}$, then since $t\in \mathfrak{m}(\eta^{-1})$ we have \begin{equation} \label{end} |t/L-a/q|\geq 1/\eta L. \end{equation}  It then follows from Lemma \ref{major}, (\ref{aM2}), (\ref{int}), and (\ref{end}) that \begin{equation*} |S(t)| \ll \eta, \end{equation*} as required.
\qed

\subsection{Proof of (\ref{moment})}
We first note that 
\begin{align*} \sum_{t \in \Z_{L}} &|S(t)|^6 = \frac{1}{M^{12}}\sum_{j\leq x_1,\dots,x_6\leq M}x_1\cdots x_6\sum_{t \in \Z_{L}}e^{2\pi i (h(x_1)+h(x_2)+h(x_3)-h(x_4)-h(x_5)-h(x_6))t/L}\\\\
&\leq \frac{L}{M^6} \cdot \#\left\{(x_1,\dots,x_6): j\leq x_i \leq M, h(x_1)+h(x_2)+h(x_3)\equiv h(x_4)+h(x_5)+h(x_6)\!\!\!\! \mod{L} \right\}.\\
\end{align*}

\noindent By definition of $j$ and $M$, both sides of the congruence above lie in $[1,L)$, so congruence modulo $L$ implies equality. Noting this fact, we have 
\begin{equation*} \sum_{t \in \Z_{L}} |S(t)|^6 \leq  \frac{L}{M^6} \cdot J(\alpha,\beta,M)
\end{equation*}
where
\begin{equation*}J(\alpha,\beta,M)=\#\{(x_1,\dots,x_6): 1\leq x_i\leq M, \alpha(x_1^2+x_2^2+x_3^2-x_4^2-x_5^2-x_6^2)=\beta(x_1+x_2+x_3-x_4-x_5-x_6)\}.
\end{equation*}
By Proposition \ref{content}, we know that $(\alpha,\beta)=\textnormal{cont}(h)\leq C$, so by (\ref{aM2}) it suffices to show under the assumption $(\alpha,\beta)=1$ that 
\begin{equation*} J(\alpha,\beta,M) \ll M^4/\alpha.
\end{equation*}
Examining the equation 
\begin{equation}\label{Jab} \alpha(x_1^2+x_2^2+x_3^2-x_4^2-x_5^2-x_6^2)=\beta(x_1+x_2+x_3-x_4-x_5-x_6),
\end{equation}
we see that the right hand side must be divisible by $\alpha$, so if $(\alpha,\beta)=1$, it must be the case that $\alpha$ divides $x_1+x_2+x_3-x_4-x_5-x_6$. Since this expression takes values in $(-3M,3M)$, there are at most $6M/\alpha +1 \leq 7M/\alpha$ choices for its value, where the last inequality follows from (\ref{cofmax}). Also, a chosen value for this expression determines the value of $x_1^2+x_2^2+x_3^2-x_4^2-x_5^2-x_6^2$ required to satisfy (\ref{Jab}).

Now we invoke a special case of the solution to Tarry's problem, which says that for any fixed $s,t \in \Z$, the number of solutions to the system 
\begin{align*}
&x_1+x_2+x_3-x_4-x_5-x_6 = s\\&  x_1^2+x_2^2+x_3^2-x_4^2-x_5^2-x_6^2=t
\end{align*}
with $1\leq x_i \leq M$ is at most $CM^3$. Discussions of this fact and Tarry's problem in general can be found in \cite{Hua} and \cite{Wooley}. Putting the pieces together, we have that for $(\alpha,\beta)=1$,
\begin{equation*}
J(\alpha,\beta,M) \ll \frac{M}{\alpha}\cdot M^3=\frac{M^4}{\alpha},
\end{equation*}
and the result follows.
\qed

\section{Proofs of Propositions \ref{intquad} and \ref{content}}

\subsection{Proof of Proposition \ref{intquad}}\label{B1}

 First we recall that a polynomial is intersective if and only if it has a root in the $p$-adic integers for every prime $p$.

\noindent Suppose $f(x)=ax^2+bx+c \in \Z[x]$ has no rational roots, hence $b^2-4ac$ is not a perfect square. Let $p=3$ if $b^2-4ac=-n^2$ for $n\in \N$, and otherwise let $p$ be any prime such that ord$_p(b^2-4ac)$, the exponent of $p$ in the prime factorization of $b^2-4ac$, is odd. Letting $\Q_p$ denote the field of $p$-adic numbers, we have that $b^2-4ac$ is not a square in $\Q_p$. Therefore, by the quadratic formula, $f$ has no roots in $\Q_p$, hence no $p$-adic integer roots, so $f$ is not an intersective polynomial.

\noindent Now suppose that $f(x)=a(\alpha x + \beta)(\gamma x + \lambda)$ with $a,\alpha,\beta,\gamma,\lambda \in \Z$ and $(\alpha,\beta)=(\gamma,\lambda)=1$. If $p$ is a prime that divides both $\alpha$ and $\gamma$, then we see that $f$ has no root modulo $p^k$ whenever $p^k \nmid a$, hence $f$ is not an intersective polynomial. 

\noindent Conversely, if $(\alpha,\gamma)=1$, we see that $-\beta /\alpha$ is a $p$-adic integer root of $f$ whenever $p \nmid \alpha$, and $-\lambda / \gamma$ is a $p$-adic integer root of $f$ whenever $p \nmid \gamma$. Since at least one of these divisibility conditions holds for every prime $p$, $f$ is an intersective polynomial.  
\qed

\subsection{Proof of Proposition \ref{content}}\label{B2}

Recall that at this stage we have a fixed intersective quadratic $f \in \Z[x]$. 

\noindent We can assume $f(x)=ax^2+bx+c=(\alpha x + \beta)(\gamma x + \lambda)$ with $(a,b,c)=(\alpha,\beta)=(\gamma,\lambda)=1$, since cont$(f_d)$ and the expression $|\alpha\lambda-\beta\gamma|$ both behave predictably under scaling of $f$. Further, since $f$ is intersective, we have that $(a,b)=(\alpha, \gamma)=1$.

\noindent In this case, we have $f_d(x)=dax^2+(2ar_d+b)x+f(r_d)/d$, so \begin{equation*}\textnormal{cont}(f_d)=(da,2ar_d+b)=(d,2ar_d+b), \end{equation*} where the last equality holds because $(a,b)=1$ implies $(a,2ar_d+b)=1$.

\noindent Now suppose that a prime power $p^k$ divides both $d$ and $2ar_d+b=f'(r_d)=\alpha(\gamma r_d+\lambda)+\gamma(\alpha r_d+\beta)$. Because $p \mid d$, and by the construction of the root $r_d$ described following Definition \ref{aux}, it is either the case that $p^k \mid \alpha r_d + \beta$ or $p^k \mid \gamma r_d + \lambda$. We will assume the former without loss of generality, so in particular $p \nmid \alpha$. 

\noindent We then see that $p^k \mid \alpha r_d + \beta$ and $p^k \mid f'(r_d)$ implies $p^k \mid \alpha (\gamma r_d + \lambda)$, and since $p \nmid \alpha$, it must be the case that $p^k \mid \gamma r_d + \lambda$. In summary, we have that \begin{equation*} r_d \equiv -\beta/\alpha \equiv -\lambda/\gamma \text{ mod }p^k. \end{equation*} In particular $p^k \mid \alpha\lambda - \beta\gamma$, and the result follows. 
\qed


\begin{thebibliography}{10}       
  
\bibitem{BPPS} 
{\sc A. Balog, J. Pelik\'an, J. Pintz, E. Szemer\'edi}, {\em Difference sets without $k$-th powers}, Acta. Math. Hungar. 65 (2) (1994), pp. 165-187.

\bibitem{BB} {\sc D. Berend, Y. Bilu}, {\em Polynomials with roots modulo every integer}, Proc. Amer. Math. Soc. 124 (1996), pp. 1663-1671.
 
\bibitem{Furst}
{\sc H.~Furstenberg}, {\em Ergodic behavior of diagonal measures and a theorem of {S}zemer\'edi on arithmetic progressions},
  J. d'Analyse Math 71 (1977), pp. 204-256.
   
\bibitem{Hua}
{\sc L. K. Hua}, {\em Additive theory of prime numbers }, American Mathematical Society, Providence, RI 1965.
  
\bibitem{KMF}
{\sc T. Kamae, M. Mend\`es France}, {\em van der Corput's difference theorem}, Israel J. Math. 31, no. 3-4, (1978), pp. 335-342.

\bibitem{Le}
{\sc T. H. L\^{e}}, {\em Intersective polynomials and the primes}, J. Number Theory 130 no. 8 (2010), pp. 1705-1717.

\bibitem{Lucier}
{\sc J. Lucier}, {\em Intersective sets given by a polynomial}, Acta Arith. 123 (2006), pp. 57-95.
 
\bibitem{LM}  
{\sc N. Lyall, \'A. Magyar}, {\em Polynomial configurations in difference sets}, J. Number Theory 129 (2009), pp. 439-450.

\bibitem{LM2}
{\sc N. Lyall, \'A. Magyar}, {\em Simultaneous polynomial
recurrence},  Bull. Lond. Math. Soc.  43  (2011),  no. 4, 765-785
 
\bibitem{Mont}
{\sc H. L. Montgomery} {\em Ten lectures on the interface between analytic number theory and harmonic analysis}, CBMS Regional Conference Series in Mathematics, 84.

\bibitem{PSS}
{\sc J. Pintz, W. L. Steiger, E. Szemer\'edi}, {\em On sets of natural numbers whose difference set contains no squares}, J. London Math. Soc. 37 (1988), pp. 219-231.

\bibitem{Roth}
{\sc K. F. Roth}, {\em On certain sets of integers}, J. London Math. Soc. 28 (1953), pp. 104-109.

\bibitem{sarkozy}
{\sc A. S\'ark\"ozy}, {\em On difference sets of sequences of integers I}, Acta. Math. Hungar. 31 (1-2) (1978), pp. 125-149.

\bibitem{Slip} {\sc S. Slijep\v{c}evi\'c}, {\em A polynomial S\'ark\"ozy-Furstenberg theorem with upper bounds}, Acta Math. Hungar. 98 (2003), pp. 275-280


\bibitem{julia}
{\sc J. Wolf}, {\em Arithmetic structures in sets of integers}, Ph.D. thesis, University of Cambridge, submitted December 2007  

\bibitem{Wooley}
{\sc T. Wooley}, {\em Some remarks on Vinogradov's mean value theorem and Tarry's problem}, Monatsh. Math. 122, no. 3 (1996), pp. 265-273. 

\end{thebibliography}
\end{document}